\documentclass{amsart}
\usepackage{amsfonts}
\usepackage{amsmath,amssymb}
\usepackage{amsthm}
\usepackage{amscd}
\usepackage{graphics}
\usepackage{graphicx}

\theoremstyle{remark}{
\newtheorem{Def}{{\rm Definition}}
\newtheorem{Ex}{{\rm Example}}
\newtheorem{Rem}{{\rm Remark}}
\newtheorem{Prob}{{\rm Problem}}

}
\theoremstyle{plain}
{

\newtheorem{Prop}{Proposition}
\newtheorem{Thm}{Theorem}
\newtheorem{MainThm}{Main Theorem}

}

\begin{document}
\title[Real algebraic functions with prescribed preimages]{Construction of real algebraic functions with prescribed preimages}
\author{Naoki kitazawa}
\keywords{Real algebraic manifolds, Real algebraic functions. Reeb graphs. Special generic maps. Morse-Bott functions. \\
\indent {\it \textup{2020} Mathematics Subject Classification}: Primary~14P05, 14P25, 57R45, 58C05. Secondary~57R19.}
\address{Institute of Mathematics for Industry, Kyushu University, 744 Motooka, Nishi-ku Fukuoka 819-0395, Japan\\
 TEL (Office): +81-92-802-4402 \\
 FAX (Office): +81-92-802-4405 \\
Osaka Central
	Advanced Mathematical Institute, 3-3-138 Sugimoto, Sumiyoshi-ku Osaka 558-8585 \\
	TEL (Office): +81-6-6605-3103 \\
	FAX (Office): +81-6-6605-3104 \\
}
\email{n-kitazawa@imi.kyushu-u.ac.jp \/ naokikitazawa.formath@gmail.com}
\urladdr{https://naokikitazawa.github.io/NaokiKitazawa.html}
\maketitle
\begin{abstract}

We present real algebraic functions with prescribed preimages.

Nash and Tognoli show that smooth closed manifolds can be the zero sets of some real polynomial maps and {\it non-singular}. The canonical projections of spheres naturally embedded in the $1$-dimensional higher Euclidean spaces and some natural functions on projective spaces, Lie groups and their quotient spaces are important examples of real algebraic functions being also Morse. In general, it is difficult to construct such examples of maps and the structures of the manifolds. In addition the maps are hard to understand globally.

We construct examples by answering to a problem from singularity theory and differential topology. It asks whether we can reconstruct nice smooth functions with prescribed preimages. We have previously given an answer with real algebraic functions. This previous result is one of our key ingredients.

\end{abstract}
\section{Introduction.}
\label{sec:1}

One of Nash and Tognoli's greatest work \cite{nash, tognoli} says that every smooth closed manifold has the structure of a so-called {\it non-singular} real algebraic manifold and the zero set of some real polynomial map. For related surveys, see \cite{kollar} for example. For some smooth closed manifolds, explicit real algebraic functions are well-known. Spheres, some projective spaces, some Lie groups and some of their quotient spaces are known to admit some natural real algebraic functions which are also Morse. The unit spheres and their canonical projections give simplest examples. In general, explicit construction of real algebraic functions and maps and knowing preimages and more generally, their explicit global structures including information on some important polynomials, are hard. We challenge the construction by applying singularity theory, differential topology, and (elementary) algebraic geometry.

\subsection{Terminologies and notation on polyhedra, manifolds, smooth maps and real algebraic maps.}

For a topological space $X$ homeomorphic to a cell complex whose dimensions is finite, we can define the dimension $\dim X$ as an integer uniquely. 
A topological space homeomorphic to a topological manifold has the structure of a CW complex. A smooth manifold is known to have the structure of a canonically and uniquely obtained polyhedron. This is a so-called PL manifold. It is also well-known that a topological space having the structure of a polyhedron of dimension at most $2$ has the structure of a polyhedron uniquely. For a topological manifold of dimension at most $3$, this also holds. For this, see \cite{moise} for example. 
  
Let ${\mathbb{R}}^k$ denote the $k$-dimensional Euclidean space. This is a smooth manifold. This is also a Riemannian manifold with the standard Euclidean metric. For a point $x \in {\mathbb{R}}^k$, we can define $||x|| \geq 0$ as the distance between $x$ and the origin $0$ under this metric. This is also a ({\it non-singular}) real algebraic manifold: the $k$-dimensional {\it real affine space}. Let $S^k:=\{x \in {\mathbb{R}}^{k+1} \mid ||x||=1\}$ denote the $k$-dimensional unit sphere. It is a smooth compact submanifold of ${\mathbb{R}}^{k+1}$ with no boundary and of dimension $k \geq 0$. It is connected for $k \geq 1$ and a discrete set with exactly two points for $k=0$. It is also a ({\it non-singular}) real algebraic submanifold and the zero set of the real polynomial ${||x||}^2-1={\Sigma}_{j=1}^{k+1} {x_j}^2 -1$ where $x:=(x_1,\cdots,x_{k+1})$.
Let $D^k:=\{x \in {\mathbb{R}}^{k} \mid ||x|| \leq 1\}$ denote the $k$-dimensional unit disk. It is a smooth compact and connected submanifold of ${\mathbb{R}}^{k}$ and of dimension $k \geq 1$.

For a differentiable map $c:X \rightarrow Y$ between differentiable manifolds, $x \in X$ is a {\it singular} point of $c$ if the rank of the differential at $x$ is smaller than the minimum in $\{\dim X, \dim Y\}$. 
For a singular point $x \in X$, $c(x)$ is a {\it singular value} of $c$.
Let $S(c)$ denote the {\it singular set} of $c$, the set of all singular points of $c$. Hereafter, we consider smooth maps or maps of
the class $C^{\infty}$ as differentiable maps unless otherwise stated. The {\it canonical projection of the unit sphere} $S^{k-1}$ is defined as the restriction of the canonical projection ${\pi}_{k,k_1}:{\mathbb{R}}^{k} \rightarrow {\mathbb{R}}^{k_1}$, mapping $x=(x_1,x_2) \in {\mathbb{R}}^{k_1} \times {\mathbb{R}}^{k_2}={\mathbb{R}}^k$ to $x_1 \in {\mathbb{R}}^{k_1}$, to $S^{k-1}$ for $k_1, k_2>0$ and $k=k_1+k_2$. This map is Morse for $k_1=1$.
As presented before, more generally, natural real algebraic Morse functions on real projective spaces, complex ones and quaternion ones are well-known. Some Lie groups are embedded naturally in Euclidean spaces and admit Morse functions represented by real polynomials. We can say this for their various quotient spaces. For related theory, see \cite{ramanujam, takeuchi} and see also \cite{maciasvirgospereirasaez} for example.
\subsection{Graphs and Reeb graphs.} 

Graphs are important tools. A graph is a $1$-dimensional CW complex with the vertex set, defined as the set of all $0$-dimensional cells and the edge set, defined as the set of all $1$-dimensional cells. A vertex and an edge are elements of these sets respectively. The closure of an edge homeomorphic to a circle is called a {\it loop}. Hereafter, a graph
has no loops and it may be a multi-graph or a graph with more than one edge connecting given two distinct vertices. 
An {\it isomorphism} from a graph $K_1$ onto another graph $K_2$ is a piecewise smooth homeomorphism mapping the edge set and the vertex set of $K_1$ onto those of $K_2$.

{\it Reeb graphs} are graphs and fundamental tools in our study. For a smooth function $c:X \rightarrow \mathbb{R}$, we can define an equivalence relation on $X$ as follows. Two points $x_1$ and $x_2$ in $X$ are equivalent if and only if they are in a same connected component of a preimage $c^{-1}(y)$ ($y \in Y$). By this equivalence relation ${\sim}_{c}$, we can define the quotient space $W_c:=X/{{\sim}_{c}}$ and the quotient map $q_c:X \rightarrow W_c$. Theorem 3.1 of \cite{saeki2} says that for a smooth function $c$
on a closed manifold with finitely many singular values, $W_c$ is a graph whose vertex set consists of all points $p$ whose preimages ${q_c}^{-1}(p)$ contain some singular points of $c$. 
{\rm Morse}({\rm -Bott}) functions and smooth functions in some considerably wide classes satisfy this. 

\begin{Def}
The graph $W_c$ is called the {\it Reeb graph}
of $c$.
\end{Def}


Reeb graphs have information on the manifolds roughly and do not miss important information. \cite{reeb} is one of pioneering studies on Reeb graphs. 
\subsection{Our problems and our main results.}
Related to our work, we present problems and history on reconstructing nice smooth functions from graphs.
\begin{Prob}
	\label{prob:1}
For a graph, can we construct a nice smooth function whose Reeb graph is isomorphic to it? We do not fix the manifold beforehand.
\end{Prob}
\cite{sharko} asked this first and smooth functions on closed surfaces have been explicitly constructed for graphs satisfying some nice conditions. \cite{masumotosaeki} generalizes this for arbitrary graphs.
\cite{martinezalfaromezasarmientooliveira, michalak} have set explicit problems and solved. Their studies are essentially on smooth functions on closed surfaces and Morse functions such that connected components of preimages containing no singular points are always spheres. The following is a revised problem, introduced first by the author in \cite{kitazawa1}.
\begin{Prob}
	\label{prob:2}
	Can we construct a smooth function whose Reeb graph is isomorphic to a given graph, whose singular points are mild in suitable senses, and whose preimages are as prescribed? We do not fix the manifold beforehand.
\end{Prob}
 \cite{kitazawa2, kitazawa6, kitazawa7} give answers (Theorems \ref{thm:3} and \ref{thm:4}). \cite{saeki2} respects some of our informal discussions on \cite{kitazawa2}.
 These are from the smooth category and  
 differential topology.
 
\begin{Prob}
\label{prob:3}
Can we construct real algebraic functions whose Reeb graphs and preimages are as prescribed as Problems \ref{prob:1} and \ref{prob:2} ask?
\end{Prob}
\cite{kitazawa3} is a pioneering answer to this, respecting Problem \ref{prob:1}, and gives functions whose Reeb graphs are isomorphic to graphs nicely embedded in ${\mathbb{R}}^n$ (Theorem \ref{thm:5}).
Preimages containing no singular points are disjoint unions of spheres for most of these functions. They generalize the canonical projections of the unit spheres.

In our study, we have the following. This constructs smooth real algebraic functions with prescribed preimages which may not be spheres and this with additional Main Theorems, presented later, respects Problem \ref{prob:2} in Problem \ref{prob:3} first. For this, see also the abstract \cite{kitazawa4} of our related
talk in a conference. Real algebraic manifolds are unions of connected components the zero sets of some real polynomial maps in our paper. {\it Non-singular} real algebraic manifolds are defined via implicit function theorem: we use the ranks of the maps defined canonically from the polynomials. 

\begin{MainThm}
	\label{mthm:1}
	Let $a>3$ and $m>2$ be integers. Let $G$ be a graph as follows.
	\begin{itemize}
		\item The vertex set is of size $l$ and the $j$-th vertex is denoted by $v_j$ for $1 \leq j \leq a$.
		\item The edge set is of size $l-1$ and the $j$-th edge connects the vertex $v_j$ and $v_{j+1}$ for $1 \leq j \leq a-1$.  
	\end{itemize}
 Let $\{F_j\}_{j=1}^{a-1}$ be a family of smooth manifolds satisfying the following conditions.
\begin{itemize}
\item The two manifolds $F_1$ and $F_{a-1}$ are the {\rm (}$m-1${\rm )}-dimensional unit spheres $S^{m-1}$. The others are the unit spheres or represented as connected sums of finitely many copies of manifolds represented as the products $S^j \times S^{m-j-1}$ with integers $1 \leq j \leq m-2${\rm :} the connected sum is taken in the smooth category.
\item For each integer $1 \leq j \leq a-2$, either $F_j$ or $F_{j+1}$ is not diffeomorphic to the unit sphere.
\end{itemize}
	Then there exist an $m$-dimensional non-singular real algebraic closed and connected manifold $M$, a smooth real algebraic function $f:M \rightarrow {\mathbb{R}}$ which is also Morse and an isomorphism $\phi:G \rightarrow W_f$ of the graphs and for the $j$-th edge $e_j$ and each point $p_{e_j}$ in the interior of the edge $\phi(e_j) \subset W_f$, the preimage ${q_f}^{-1}(p_{e_j})$ is diffeomorphic to $F_j$. $F_j$ is regarded as a non-singular real algebraic hypersurface of $M$.
	\end{MainThm}


In the next section, we prove Main Theorems.

\noindent {\bf Conflict of interest.} \\
The author was a member of the project JSPS Grant Number JP17H06128. He was also a member of the project JSPS KAKENHI Grant Number JP22K18267 "Visualizing twists in data through monodromy" (Principal Investigator: Osamu Saeki). Our present study thanks the project for their support. 
The author works at Institute of Mathematics for Industry (https://www.jgmi.kyushu-u.ac.jp/en/about/young-mentors/). Our study thanks this project.
The author is a researcher at Osaka Central
Advanced Mathematical Institute (OCAMI researcher) whereas he is not employed there. Our study also thanks
this for giving us the opportunity to study further. \\ \\
\ \\
{\bf Data availability.} \\
Data essentially supporting our present study are all contained in the present paper.
\section{On Main Theorems.}
\subsection{Additional several terminologies, notions and notation.}
A {\it diffeomorphism} means a smooth homeomorphism with no singular points. A {\it diffeomorphism on a smooth manifold} means a diffeomorphism from it to itself.
The {\it diffeomorphism type} of a smooth manifold is defined as the equivalence class under the natural equivalence relation on the family of all smooth manifolds defined by the existence of diffeomorphisms.

The {\it diffeomorphism group} of a smooth manifold is the group of all diffeomorphisms on it. This is also a topological group and topologized with the so-called {\it Whitney $C^{\infty}$ topology}. More generally, {\it Whitney $C^{\infty}$ topologies} on the set of all smooth maps between two given smooth manifolds and subspaces of this space are important in the (singularity) theory of smooth maps for example. 
For singularity theory, refer to \cite{golubitskyguillemin} for example.

A {\it smooth} bundle means a bundle whose fiber is a smooth manifold and whose structure group is regarded as (some subgroup) of the diffeomorphism group of the fiber.

We introduce {\it fold} maps.

\begin{Def}
	Let $X$ and $Y$ be smooth manifolds with no boundaries satisfying $\dim X \geq \dim Y$.
	A {\it fold} map $c:X \rightarrow Y$ is a smooth map such that at each singular point $p$, we have a suitable integer $0 \leq i(p) \leq \frac{\dim X-\dim Y+1}{2}$, local coordinates around $p$ and $c(p)$, and a local form \\
$c(x_1,\cdots,x_{\dim X})=(x_1,\cdots,x_{\dim Y-1},{\Sigma}_{j=1}^{\dim X-\dim Y-i(p)+1} {x_{\dim Y-1+j}}^2-{\Sigma}_{j=1}^{i(p)} {x_{\dim X-i(p)+j}}^2)$.
\end{Def}
\begin{Prop}
	In the previous definition, $i(p)$ is chosen uniquely and defined as the {\rm index} of $p$. The set of all singular points of $c$ of a fixed index is a smooth regular submanifold of $c$, with no boundary, and of dimension $\dim Y-1$. If $X$ is closed, then the submanifold is compact. The restriction to the submanifold is a smooth immersion. 
\end{Prop}
\begin{Def}
	If in the definition of a fold map, $i(p)=0$ always holds, then this is called a {\it special generic} map.
\end{Def}
A Morse function is of course a fold map. 
In short, fold maps are locally projections or the product map of a Morse function and the identity map on some disk. For special generic maps, this local Morse function is chosen as a so-called {\it height function} of a unit disk. A {\it height function} $h$ of a unit disk is a function of the form $h(x)=\pm ||x||^2+c$ where $c$ is some real number. 
For fold maps, see also \cite{golubitskyguillemin} for example.
We introduce very fundamental and explicit special generic maps, discussed in \cite{saeki1} for example. They are also key tools in our main result.

\begin{Ex}
\label{ex:1}
	\begin{enumerate}
		\item \label{ex:1.1}
		 The canonical projections of the unit spheres are special generic. The restrictions to the singular sets, which are also regarded as the unit spheres, are embeddings. The images are regarded as the unit disks whose dimensions are same as those of the Euclidean spaces of the targets.
		\item \label{ex:1.2}
		 Let $m \geq n \geq 2$ be integers. Let $M$ be an $m$-dimensional smooth manifold diffeomorphic to one represented as a connected sum of $l> 0$ manifolds diffeomorphic to $S^{k_j} \times S^{m-k_j}$ for each integer $1 \leq j \leq l$ and some integer $1 \leq k_j \leq n-1$ where the connected sum is taken in the smooth category. We easily have a special generic map $f:M \rightarrow {\mathbb{R}}^n$ such that the restriction to the singular set $S(f)$ is an embedding and that the image is a smoothly embedded submanifold diffeomorphic to one represented as a boundary connected sum of $l> 0$ manifolds diffeomorphic to $S^{k_j} \times D^{n-k_j}$ for each integer $1 \leq j \leq l$. The boundary connected sum is, as before, taken in the smooth category. 
		\end{enumerate}

\end{Ex}
We can know the construction of these special generic maps from fundamental arguments of \cite{saeki1}.
These special generic maps are generalized to the class of Definition \ref{def:4}.
This class of special generic maps is first defined in \cite{kitazawa9} and later renamed in \cite{kitazawa5}.
\begin{Def}
\label{def:4}
Let $m \geq n \geq 1$ be integers. Let $M$ be an $m$-dimensional closed and connected manifold. Let $f:M \rightarrow {\mathbb{R}}^n$ be a special generic map. We assume that the restriction $f {\mid}_{S(f)}:S(f) \rightarrow {\mathbb{R}}^n$ is an embedding and that the image $D_M:=f(M)$ is an $n$-dimensional smoothly embedded manifold. We also assume the following.
\begin{itemize}
	\item There exists a small collar neighborhood $N(\partial D_M)$ of the boundary $\partial D_M$ of the image $D_M$ and the composition of the restriction of the map to the preimage $f^{-1}(N(\partial D_M))$ with the canonical projection to the boundary gives a trivial smooth bundle whose fiber is diffeomorphic to a unit disk $D^{m-n+1}$. Let $M_{f,{\rm B}}$ denote this bundle.
	\item On the complementary set $D_M-{\rm Int}\ N(\partial D_M)$ of the interior of the collar neighborhood before, the restriction $f {\mid}_{f^{-1}(D_M-{\rm Int}\ N(\partial D_M))}:f^{-1}(D_M-{\rm Int}\ N(\partial D_M)) \rightarrow D_M-{\rm Int}\ N(\partial D_M)$ of the map gives a trivial smooth bundle whose fiber is diffeomorphic to the unit sphere $S^{m-n}$. Let $M_{f,{\rm I}}$ denote this bundle. 
\item The two bundles $\partial M_{f,{\rm B}}$ and $\partial M_{f,{\rm I}}$ defined canonically on the boundaries of the previous bundles are glued by the product map of the diffeomorphism for the natural identification between the base spaces and the identity map on the fiber $S^{m-n}=\partial D^{m-n+1}$. Note also that the fibers are identified in a canonical way.
\end{itemize}
Then the map $f$ is called a {\it product-organized} special generic map.
\end{Def}
We present an explicit fact as Example \ref{ex:2}. We can also know this from fundamental arguments of \cite{saeki1}.
\begin{Ex}
\label{ex:2}
Let the image $f(M)=D_M$ of a product-organized special generic map $f:M \rightarrow {\mathbb{R}}^n$ be as presented in Example \ref{ex:1}. Then the manifold $M$ of the domain is diffeomorphic to the presented $m$-dimensional manifold in Example \ref{ex:1}.
\end{Ex}

\subsection{Proofs of Main Theorems.}


\cite{kitazawa3, kitazawa8} present important arguments and facts in our proofs. We review some of these results.



For $s_1,s_2 \in \mathbb{R}$ satisfying the relation $s_1 \leq s_2$, define $[s_1,s_2]:=\{x \mid s_1 \leq x \leq s_2\} \subset \mathbb{R}$ in a standard way.

\begin{Thm}[\cite{kitazawa3, kitazawa8}]
\label{thm:1}	
Let $D \subset {\mathbb{R}}^n$ be a connected open subset in ${\mathbb{R}}^n$ such that the closure $\overline{D}$ is compact and that the boundary $\partial \overline{D}$ of the closure consists of exactly $l \geq 1$ mutually disjoint non-singular real 
algebraic hypersurfaces {\rm (}{\rm (}$n-1${\rm )}-dimensional non-singular real algebraic manifolds in ${\mathbb{R}}^n${\rm )}. Let $\{S_j\}_{j=1}^l$ denote the family of these $l$ hypersurfaces. Suppose the following conditions.
\begin{itemize}
\item The set $S_j$ is a connected component of the zero set of some real polynomial $f_j(x_1, \cdots x_n)$.
\item There exists a connected open neighborhood $U$ of $\overline{D}$ and for each integer $1 \leq j \leq l$, except $S_j$, connected components of the zero set of the real polynomial $f_j(x_1, \cdots x_n)$ is disjoint from $U$. 
\item It holds that $D \subset {\bigcap}_{j=1}^l \{(x_1, \cdots x_n) \in {\mathbb{R}}^n \mid f_j(x_1, \cdots x_n)>0\}$. Two relations $U \bigcap {\bigcap}_{j=1}^l \{(x_1, \cdots x_n)  \in {\mathbb{R}}^n \mid f_j(x_1, \cdots x_n)>0\}=D$ and $U \bigcap {\bigcap}_{j=1}^l \{(x_1, \cdots x_n) \in {\mathbb{R}}^n \mid f_j(x_1, \cdots x_n) \geq 0\}=\overline{D}$ also hold. 
\end{itemize}
Then for any integer $n^{\prime}>n$, there exists another connected open subset $D^{\prime} \subset {\mathbb{R}}^{n^{\prime}}$ in ${\mathbb{R}}^{n^{\prime}}$ such that the closure $\overline{D^{\prime}}$ is compact and that the boundary $\partial \overline{D^{\prime}}$ of the closure is a non-singular real algebraic connected hypersurface enjoying the following properties.
\begin{enumerate}
\item
\label{thm:1.1}

The closure $\overline{D}$ of $D$ is embedded in $\overline{D^{\prime}}$ by the inclusion mapping $x \in \overline{D} \subset {\mathbb{R}}^n$ to $(x,0) \in \overline{D^{\prime}} \subset {\mathbb{R}}^{n^{\prime}}= {\mathbb{R}}^n \times  {\mathbb{R}}^{n^{\prime}-n}$. The inclusion maps the boundary $\partial \overline{D} \subset \overline{D}$ into the boundary $\partial \overline{D^{\prime}} \subset \overline{D^{\prime}}$ and the interior $D$ into the interior $D^{\prime}$.
\item
\label{thm:1.2}
The boundary $\partial \overline{{D^{\prime}}}$ is a connected component of the zero set of the real polynomial $F(x_1,\cdots x_n, \cdots x_{n^{\prime}}):={\prod}_{j=1}^l (f_j(x_1,\cdots x_n))-{
\Sigma}_{j=1}^{n^{\prime}-n} {x_{n+j}}^2$ and it also holds that $\partial \overline{{D^{\prime}}}=\{(x_1, \cdots x_n) \in \overline{D} \mid F(x_1,\cdots x_n, \cdots x_{n^{\prime}})=0\}$. It holds that $D^{\prime}=\{(x_1, \cdots x_n) \in D \mid F(x_1,\cdots x_n, \cdots x_{n^{\prime}})>0\}=\{(x_1, \cdots x_n) \in U \mid F(x_1,\cdots x_n, \cdots x_{n^{\prime}})>0\} \subset \{(x_1, \cdots x_n)  \in {\mathbb{R}}^n \mid F(x_1,\cdots x_n, \cdots x_{n^{\prime}})>0\}$. The connected open set $D^{\prime}$ is the bounded connected component of the set ${\mathbb{R}}^{n^{\prime}}-\partial \overline{{D^{\prime}}}$ and uniquely determined. There exists a connected open neighborhood $U^{\prime}$ of $\overline{D^{\prime}}$ and two relations $U^{\prime} \bigcap \{(x_1,\cdots x_n, \cdots x_{n^{\prime}}) \in {\mathbb{R}}^{n^{\prime}}  \mid F(x_1, \cdots x_n, \cdots x_{n^{\prime}})>0\}=D^{\prime}$ and $U^{\prime} \bigcap \{(x_1,\cdots x_n, \cdots x_{n^{\prime}}) \in {\mathbb{R}}^{n^{\prime}} \mid F(x_1,\cdots x_n, \cdots x_{n^{\prime}}) \geq 0\}=\overline{D^{\prime}}$ hold. 
\item
\label{thm:1.3}
 The restriction ${\pi}_{n^{\prime},n} {\mid}_{\partial \overline{{D^{\prime}}}}:\partial \overline{{D^{\prime}}} \rightarrow {\mathbb{R}}^n$ of the canonical projection ${\pi}_{n^{\prime},n}$ to the boundary $\partial \overline{{D^{\prime}}}$ is a product-organized special generic map.
\item
\label{thm:1.4}
The smooth manifold $\overline{D^{\prime}}$ is diffeomorphic to the manifold obtained in the following steps.
\begin{enumerate}
\item
\label{thm:1.4.1}
Prepare a smooth manifold $P_1:=\overline{D^{\prime}} \times D^{n^{\prime}-n}$. We regard this as a trivial smooth bundle over $\overline{D^{\prime}}$.
\item
\label{thm:1.4.2}
Prepare a smooth manifold $P_2:=(\partial \overline{D^{\prime}}) \times D^{n^{\prime}-n+1}$ and choose a copy ${D_0}^{n^{\prime}-n}$ of the disk $D^{n^{\prime}-n}$ smoothly embedded in the boundary $\partial D^{n^{\prime}-n+1}$ of the disk $D^{n^{\prime}-n+1}$. We regard this as a trivial smooth bundle over the boundary $\partial \overline{D^{\prime}}$ of $\overline{D^{\prime}}$.
\item
\label{thm:1.4.3}
Glue $P_1$ and $P_2$ by a bundle isomorphism between the following trivial bundles $B_1$ and $B_2$ represented as the product map of the identity map on the base space and the identity map on the fiber. The former bundle $B_1$ is the restriction of the trivial bundle $P_1$ to the boundary $\partial \overline{D^{\prime}}$ of the base space $\overline{D^{\prime}}$. The latter bundle $B_2$ is a subbundle of the trivial bundle $P_2$ whose fiber is the subspace ${D_0}^{n^{\prime}-n} \subset \partial D^{n^{\prime}-n+1}$. Note that the fibers $D^{n^{\prime}-n}$ and ${D_0}^{n^{\prime}-n}$ are identified in a canonical way.
\item 
\label{thm:1.4.4}
Eliminate the corner of the resulting manifold.
\end{enumerate}
\item
\label{thm:1.5}
In the case $U={\mathbb{R}}^n$, we can also put $U^{\prime}={\mathbb{R}}^{n^{\prime}}$. Furthermore, the set $\partial \overline{D^{\prime}} \subset {\mathbb{R}}^{n^{\prime}}$ is also the zero set of the real polynomial $F$.
\end{enumerate}

\end{Thm}
We review main ingredients of the original proof of this.
\begin{proof}[Reviewing main ingredients of the original proof of Theorem \ref{thm:1}.]

First, It is discussed that the open connected set $D^{\prime} \subset {\mathbb{R}}^{n^{\prime}}$ is a connected open set the boundary $\partial \overline{D^{\prime}}$ of whose closure $\overline{D^{\prime}}$ is an ($n^{\prime}-1$)-dimensional non-singular real algebraic closed and connected manifold in \cite{kitazawa3}. The properties (\ref{thm:1.1}) and (\ref{thm:1.2}) are also discussed in \cite{kitazawa3} except the statement on the set $U^{\prime}$. In addition the definition of the set $U^{\prime}$ is not explicitly discussed in \cite{kitazawa8}. Consult the paper \cite{kitazawa3} mainly for the presented properties. Note that the connected open set $U$ is chosen as the whole set ${\mathbb{R}}^n$ in \cite{kitazawa3}. We can replace it by a general connected open set $U \subset {\mathbb{R}}^n$ being also a neighborhood of the set $\overline{D} \subset {\mathbb{R}}^n$.
 
By our construction, we can put $U^{\prime}=U \times {\mathbb{R}}^{n^{\prime}-n} \subset {\mathbb{R}}^{n^{\prime}}$ for example. From this argument, we also have the property (\ref{thm:1.5}).

 Most of our remaining presentation on the proof is presented first in the preprint \cite{kitazawa8}. For our understanding, we present related arguments. We can understand these arguments in self-contained ways except "Discussion 14 in \cite{kollar}". 

For a sufficiently small positive number $a_0$ and $0 \leq t \leq 1$, we can define $S_{f_t}:=U \bigcap \{x \in {\mathbb{R}}^{n} \mid {\prod}_{j=1}^l (f_j(x))-a_0t=0\} \subset \overline{D} \subset U$. The union
${\bigcup}_{t \in [0,1]} S_{f_t}$ can be naturally regarded as a small collar neighborhood of $\partial \overline{D}=S_{f_0} \subset \overline{D} \subset U$. For this, see also Discussion 14 in \cite{kollar} for example.
We can also define the set $S_{F_t}:=\{(x,x^{\prime}) \in S_{f_t} \times {\mathbb{R}}^{n^{\prime}-n} \subset {\mathbb{R}}^{n} \times {\mathbb{R}}^{n^{\prime}-n}={\mathbb{R}}^{n^{\prime}} \mid F(x,x^{\prime})=0\} \subset \partial \overline{D^{\prime}}$. We also use the notation $x=(x_1,\cdots x_n)$ and $x^{\prime}=(x_{n+1},\cdots x_{n^{\prime}})$ here. 
We study the restriction of the projection ${\pi}_{n^{\prime},n}$ to the set $\partial \overline{D^{\prime}}$. The closure $\overline{D_0}$ of the complementary set $D_0:=\overline{D}-{\bigcup}_{t \in [0,1]} S_{f_t}$ is also important. \\ 
\ \\
PART 1 The restriction of the projection ${\pi}_{n^{\prime},n}$ to the set $(\partial \overline{D^{\prime}}) \bigcap {\bigcup}_{t \in [0,1]} S_{F_t}={\bigcup}_{t \in [0,1]} S_{F_t}$.  \\
This is regarded as the product map of a height function of the unit disk $D^{n^{\prime}-n}$ and the identity map on $\partial \overline{D}$. \\
\ \\
PART 2 The restriction of the projection ${\pi}_{n^{\prime},n}$ to the set ${{\pi}_{n^{\prime},n}}^{-1}(\overline{D_0}) \bigcap \partial \overline{D^{\prime}}$. \\
The domain of the map is also same as the closure of the complementary set of ${\bigcup}_{t \in [0,1]} S_{F_t} \subset \partial \overline{D^{\prime}}$. 
As we can also see in \cite{kitazawa3}, the preimage of the closure $\overline{D_0}$ is locally regarded as the graph of a smooth map on an ($n^{\prime}-1$)-dimensional smooth manifold. Here we regard some variable $x_j$ ($n<j \leq n^{\prime}$) as the variable of the target of the function. 

By mapping $(x,x^{\prime}) \in (\overline{D_0} \times {\mathbb{R}}^{n^{\prime}-n}) \bigcap \partial \overline{D^{\prime}} \subset \partial \overline{D^{\prime}} \subset  {\mathbb{R}}^{n} \times {\mathbb{R}}^{n^{\prime}-n}={\mathbb{R}}^{n^{\prime}}$ to $(x,\frac{1}{||x^{\prime}||} x^{\prime})$, we can see that the restriction of the projection ${\pi}_{n^{\prime},n}$ gives a trivial bundle over the closure $\overline{D_0}$ of the complementary set $D_0=\overline{D}-{\bigcup}_{t \in [0,1]} S_{f_t}$ whose fiber is $S^{n^{\prime}-n-1}$.
Remember the definitions and note that $x^{\prime}$ is not the origin $0 \in{\mathbb{R}}^{n^{\prime}-n}$ for $(x,x^{\prime}) \in \overline{D_0} \times {\mathbb{R}}^{n^{\prime}-n}$. \\
\ \\
PART 1 and PART 2 show that the restriction of the canonical projection ${\pi}_{n^{\prime},n}$ to the boundary $\partial \overline{{D^{\prime}}}$ is a product-organized special generic map. We also have the property (\ref{thm:1.3}).
By our definitions, we can also check the property (\ref{thm:1.4}).
\end{proof}
In the present paper, on Theorem \ref{thm:1}, we concentrate on cases satisfying the assumption of (\ref{thm:1.5}).
These cases generalize cases where the sets $D \subset {\mathbb{R}}^n$ are surrounded by spheres $S_j$ centered at some points in ${\mathbb{R}}^n$.
\begin{Thm}
\label{thm:2}	
In Theorem \ref{thm:1}, we consider the restriction of the map ${\pi}_{n,1}$ to the boundary $\partial \overline{D}$. Suppose that it is a Morse function and that its singular set is finite and represented as $\{p_1,\cdots, p_i\} \subset \partial \overline{D}$. Then the restriction of the map ${\pi}_{n,1} \circ {\pi}_{n^{\prime},n}={\pi}_{n^{\prime},1}$ to the boundary $\partial \overline{D^{\prime}}$ is also a Morse function and its singular set is finite and equal to the preimage ${{\pi}_{n^{\prime},n}}^{-1}(\{p_1,\cdots, p_i\}) \subset {{\pi}_{n^{\prime},n}}^{-1}(\partial \overline{D})$. Furthermore, the restriction of the projection ${\pi}_{n^{\prime},n}$ to the preimage ${{\pi}_{n^{\prime},n}}^{-1}(\{p_1,\cdots, p_i\})$ gives a diffeomorphism onto the finite set $\{p_1,\cdots, p_i\}$.
\end{Thm}
We can see this fact by fundamental arguments on singularity theory of smooth maps. More precisely, we consider local coordinates around singular points of the maps and the functions.

\begin{figure}
	
	\includegraphics[height=75mm, width=100mm]{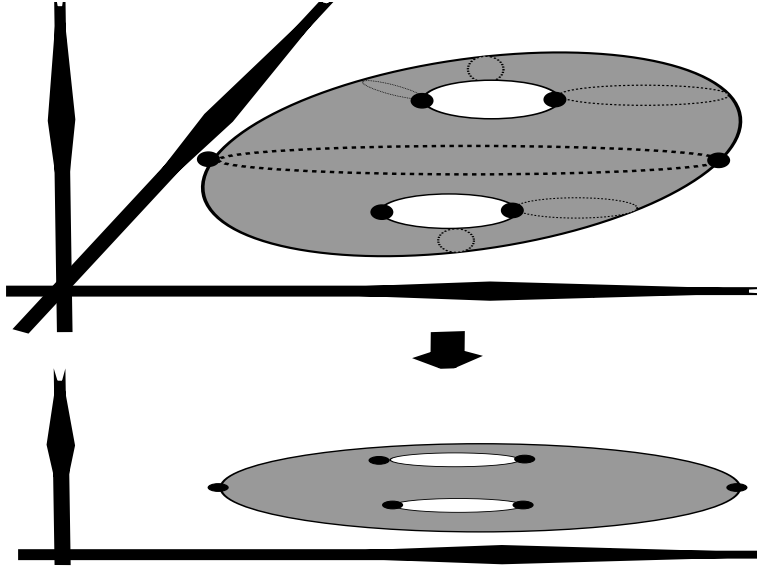}

	\caption{Example \ref{ex:3} (\ref{ex:3.2}) ($n^{\prime}=3$ and $l=3$). The upper picture is for $D^{\prime} \subset {\mathbb{R}}^{n^{\prime}}={\mathbb{R}}^3$ and the lower picture is for $D \subset {\mathbb{R}}^n={\mathbb{R}}^2$.
Dots are for singular points and singular values of the Morse function.}
	\label{fig:1}
\end{figure}

We present simplest examples on Theorems \ref{thm:1} and \ref{thm:2}. They are also from Examples \ref{ex:1} and \ref{ex:2}.
\begin{Ex}
\label{ex:3}
\begin{enumerate}
\item
\label{ex:3.1}
In Theorem \ref{thm:1}, let $D:={\rm Int}\ D^n$. We can regard $l=1$ and $S_1:=S^{n-1}$. We have $D^{\prime}={\rm Int}\ D^{n^{\prime}}$ in this case. The resulting product-organized special generic map is the canonical projection of the unit sphere $S^{n^{\prime}-1}$. By considering Theorem \ref{thm:2}, we also have a Morse function with exactly two singular points where $n \geq 2$ and this is also the canonical projection of $S^{n^{\prime}-1}$.
\item
\label{ex:3.2}
 In Theorem \ref{thm:1}, let $n=2$ and $D \subset {\mathbb{R}}^2$ a region obtained in the following way.
\begin{itemize}
\item We prepare (the interior of) a disk $D_0$ of a sufficiently large radius in ${\mathbb{R}}^2$ centered at the origin $0 \in {\mathbb{R}}^2$. 
\item We prepare $l-1>0$ disks centered at points in $D_0$ and contained in $D_0$. We also choose them in such a way that the closures are mutually disjoint.
\end{itemize}
We can regard that each $S_j$ is a circle. Then $D^{\prime} \subset {\mathbb{R}}^{n^{\prime}}$ is the interior of a smoothly embedded compact manifold in ${\mathbb{R}}^{n^{\prime}}$ diffeomorphic to one represented as a boundary connected sum of $l-1$ copies of $S^1 \times D^{n^{\prime}-1}$.

The resulting product-organized special generic map is a map on a manifold diffeomorphic to one represented as a connected sum of $l-1$ copies of $S^1 \times S^{n^{\prime}-2}$. By considering Theorem \ref{thm:2}, we also have a Morse function with exactly $2l$ singular points.
\item
\label{ex:3.3}
 This is a key ingredient in our proof of Main Theorem \ref{mthm:1}. 
This also generalizes the previous case.
In Theorem \ref{thm:1}, let $D \subset {\mathbb{R}}^n$ be represented as the interior of a smoothly embedded compact manifold in ${\mathbb{R}}^n$ diffeomorphic to one represented as a boundary connected sum of $l_0>0$ manifolds in the family $\{S^{k_j} \times D^{n-k_j}\}_{j=1}^{l_0}$ with the condition $1 \leq k_j \leq n-1$. Then $D^{\prime}$ is the interior of a smoothly embedded compact manifold in ${\mathbb{R}}^{n^{\prime}}$ diffeomorphic to one represented as a boundary connected sum of $l_0$ manifolds in the family $\{S^{k_j} \times D^{n^{\prime}-k_j}\}_{j=1}^{l_0}$.
The resulting product-organized special generic map is a map on a manifold diffeomorphic to one represented as a connected sum of $l_0$ manifolds in the family $\{S^{k_j} \times S^{n^{\prime}-k_j-1}\}_{j=1}^{l_0}$.

\end{enumerate}
\end{Ex}
We present a case of Example \ref{ex:3} (\ref{ex:3.2}) with $n^{\prime}=3$ and $l=3$ in FIGURE \ref{fig:1}.
Of course, in Example \ref{ex:3}, connected sums and boundary connected sums are taken in the smooth category.

We prove Main Theorem \ref{mthm:1}. Here an {\it ellipsoid} in ${\mathbb{R}}^k$ {\it centered at} $x_0=(x_{0,1},\cdots x_{0,k}) \in {\mathbb{R}}^k$ means a $k$-dimensional smooth compact submanifold represented as $\{(x_1,\cdots x_k) \in {\mathbb{R}}^k \mid {\Sigma}_{j=1}^k \frac{(x_j-x_{0,j})^2}{r_j}-1 \leq 0\}$ where each $r_j$ is a positive real number. This is of course diffeomorphic to the unit disk $D^k$. The unit disk $D^k$ is a specific case of such sets.

\begin{proof}[A proof of Main Theorem \ref{mthm:1}]
	
	First we prepare a connected open set $D_0$ in ${\mathbb{R}}^3$ whose closure is a disk centered at the origin $0 \in {\mathbb{R}}^3$. We also choose the disk whose radius is $R>0$. We can choose $R$ as a sufficiently large number.
	
	 We consider an increasing sequence $\{t_j\}_{j=1}^{a} \subset [-R,R]$ of real numbers of length $a$ satisfying $t_1:=-R$ and $t_{l}:=R$.
Let $t_{0,j}:=\frac{t_{j}+t_{j+1}}{2}$ for each integer $2 \leq j \leq a-2$.

       Suppose that $F_j$ is a manifold diffeomorphic to one represented as a connected sum of manifolds in the family containing exactly $k_{j,j^{\prime}} \geq 0$ copies of $S^{j^{\prime}} \times S^{m-j^{\prime}-1}$ for each integer $1 \leq j^{\prime} \leq \frac{m-1}{2}$. We also regard that $k_{j,j^{\prime}}=0$ always holds in the case $F_j$ is the unit sphere. $F_1$ and $F_{a-1}$ are the unit spheres for example. 

We consider the following procedure for each integer $0 \leq i \leq m-3$ inductively. We inductively define a connected open subset $D_i \subset {\mathbb{R}}^n$.
We start by putting $i=0$.  
\begin{itemize}
\item We check whether $i+1 \leq \frac{m-1}{2}$ or not.
\item We do the following if and only if $i+1 \leq \frac{m-1}{2}$ holds.
For each integer $2 \leq j \leq a-2$, we can choose exactly $k_{j,i+1}$ mutually disjoint connected open sets satisfying the following.

\begin{itemize}
\item The closures of the connected open sets are ellipsoids centered at points of the form $(t_{0,j},t_{1,j,j^{\prime}},0 \cdots)$ and subsets of $D_i$. For two distinct points $(t_{0,j_1},t_{1,j_1,j_{2,1}},0 \cdots)$ and $(t_{0,j_2},t_{1,j_2,j_{2,2}},0 \cdots)$ here, two numbers $t_{1,j_1,j_{2,1}}$ and $t_{1,j_2,j_{2,2}}$ are distinct. 
\item For each ellipsoid before, "$r_1$ in the definition of an ellipsoid" is $({\frac{t_{j+1}-t_{j}}{2}})^2$. 
\item By choosing the number $R>0$ as a sufficiently large one beforehand for example, we can also choose these connected open sets in such a way that the closures are mutually disjoint and we do. For this, we also need to choose each positive number $r_{j^{\prime}}$ as a sufficiently small one for $j^{\prime} \geq 2$ in the definition of an ellipsoid.
\end{itemize}

We remove the closures of the connected open sets before from $D_i$.
\item Let the resulting connected open set denoted by $D_{i,0}$. We do not change the connected open set in the case we skip the previous step and in this case we also use $D_{i,0}$.

We can apply Theorem \ref{thm:1} by putting $D:=D_{i,0}$, $n:=\dim D_i=\dim D_{i,0}$ and $n^{\prime}:=n+1$. We can also define $D_{i+1}:=D^{\prime}$ where "$D^{\prime}$" is from Theorem \ref{thm:1}. We can also apply Theorem \ref{thm:1} (\ref{thm:1.5}). We apply Theorem \ref{thm:1} in this way.
\item If $i<m-3$, then we replace the integer "$i$" by $i+1$ and go to the first step again. If for the integer "$i$", $i=m-3$, then we finish the procedure.
\end{itemize}
After this procedure, we consider the restriction of the canonical projection ${\pi}_{m+1,m}$ to the boundary $M:=\partial \overline{D^{\prime}}=\partial \overline{D_{m-2}}$, which is an $m$-dimensional non-singular real algebraic closed and connected manifold, in Theorem \ref{thm:1}. The set $\overline{D_{m-2}}$ is the closure of the open set $D_{m-2} \subset {\mathbb{R}}^{m+1}$ of course.
By composing the canonical projection ${\pi}_{m,1}$, we have a smooth real algebraic function $f:M \rightarrow \mathbb{R}$. From Theorem \ref{thm:2}, the function is also Morse. The set of all singular values of the function is $\{t_j\}_{j=1}^{a} \subset [-R,R]$.

We discuss each preimage $f^{-1}(t_{0,j})$ for $f$. This is an ($m-1$)-dimensional smooth closed and connected manifold and a non-singular real algebraic hypersurface of $M$. For each regular value $t_{\rm r}$ of $f$, $f^{-1}(t_{\rm r})$ is also a non-singular real algebraic hypersurface of the preimage ${{\pi}_{m+1,1}}^{-1}(t_{\rm r})$, regarded as a copy of the $m$-dimensional real affine space ${\mathbb{R}}^m$ embedded naturally in the ($m+1$)-dimensional real affine space ${\mathbb{R}}^{m+1}$.

To investigate the preimage $f^{-1}(t_{0,j}) \subset {{\pi}_{m+1,1}}^{-1}(t_{0,j})$, we can apply Theorem \ref{thm:1} inductively.
More precisely, we can do the following procedure for the ($i^{\prime}+1$)-th step ($0 \leq i^{\prime} \leq m-3$). 
Before starting the procedure, let $D_{t_{0,j},0}:={\rm Int}\ D^2 \subset D^2 \subset {\mathbb{R}}^2$ and put $i^{\prime}=0$. 

\begin{itemize}
\item We check whether $i^{\prime}+1 \leq \frac{m-1}{2}$ or not.
\item We do the following if and only if $i^{\prime}+1 \leq \frac{m-1}{2}$ holds.

We can choose exactly $k_{j,i^{\prime}+1}$ mutually disjoint connected open sets satisfying the following.

\begin{itemize}
\item The closures of the connected open are ellipsoids and subsets of $D_{t_{0,j},i^{\prime}}$. 
\item The closures are mutually disjoint. We choose sufficiently small ellipsoids to satisfy this condition. 
\end{itemize}
We remove the ellipsoids. 
\item Let the resulting connected open set in the underlying real affine space denoted by $D_{t_{0,j},i^{\prime},0}$. We do not change the connected open set in the case we skip the previous step and in this case we also use $D_{t_{0,j},i^{\prime},0}$.
We can apply Theorem \ref{thm:1} by putting $D$ as the resulting connected open set $D_{t_{0,j},i^{\prime},0}$, $n:=\dim D_{t_{0,j},i^{\prime}}=\dim D_{t_{0,j},i^{\prime},0}$ and $n^{\prime}:=n+1$. We can also define $D_{t_{0,j},i^{\prime}+1}:=D^{\prime} \subset {\mathbb{R}}^{n^{\prime}}={\mathbb{R}}^{n+1}$ where the set $D^{\prime}$ is abused from Theorem \ref{thm:1}. We can also apply Theorem \ref{thm:1} (\ref{thm:1.5}). We apply the theorem. 
\item If $i^{\prime}<m-3$, then we replace the integer "$i^{\prime}$" by $i^{\prime}+1$ and go to the first step again. If $i^{\prime}=m-3$, then we finish the procedure and we can see that $f^{-1}(t_{0,j})$ is diffeomorphic to the boundary $\partial \overline{D^{\prime}}=\partial \overline{D_{t_{0,j},i^{\prime}+1}} \subset \overline{D_{t_{0,j},i^{\prime}+1}}$. Of course the set $\overline{D_{t_{0,j},i^{\prime}+1}}$ is the closure of the open set $D_{t_{0,j},i^{\prime}+1}$ in (a copy of) the real affine space ${\mathbb{R}}^{m}$.
\end{itemize}

For applying Theorem \ref{thm:1} here, we can also say that here we also apply Example \ref{ex:3} (Example \ref{ex:3} (\ref{ex:3.3}) mainly). 
In the case $F_j$ is the unit sphere, the preimage $f^{-1}(t_{0,j})$ is also diffeomorphic to the unit sphere and $F_j$. This also respects $F_1$ and $F_{l-1}$ in addition to $F_j$ ($2 \leq j \leq a-2$): of course we can define the numbers $t_{0,1}$ and $t_{0,l-1}$ and the preimages in the same way.
In the case $F_j$ is not a sphere, the preimage $f^{-1}(t_{0,j})$ is also diffeomorphic to $F_j$.




We can naturally have a suitable isomorphism $\phi:G \rightarrow W_f$ and we can also see that the manifolds of the preimages are as desired.


This completes the proof.  
\end{proof}

"FIGURE 1 of \cite{kitazawa3}" shows two explicit cases for the paper.
More precisely, these two graphs show two simplest examples of {\it Poincar\'e-Reeb graphs} of {\it algebraic domains}. 
\cite{bodinpopescupampusorea} studies {\it algebraic domains} collapsing to {\it Poincar\'e-Reeb graphs} of them. An {\it algebraic domain} means an open set in a real affine space the boundary of whose closure is surrounded by mutually disjoint non-singular real algebraic hypersurfaces. A {\it Poincar\'e-Reeb} graph of an algebraic domain is a graph which the algebraic domain naturally collapses to. Rigorously, it is defined for a pair of an algebraic domain in a real affine space and the canonical projection of the real affine space to the $1$-dimensional real affine space $\mathbb{R}$.
	
We discuss the lower figure of "FIGURE 1 of \cite{kitazawa3}" and present another result. The {\it degree} of a vertex of a graph means the number of edges containing it.

\begin{MainThm}
\label{mthm:2}
Let $m>1$ be an integer. 

Let $G$ be a graph as follows. This is the lower figure of "FIGURE 1 of \cite{kitazawa3}". 
Let $b>1$ be an integer. This graph is a graph with exactly $2$ vertices of degree $1$, exactly $2$ vertices of degree $b+1$ and exactly $b+2$ edges.

\begin{itemize}
\item The first two vertices are denoted by $v_{\rm l}$ and $v_{\rm r}$, respectively. 
\item The other two vertices are denoted by $v_{\rm 1}$ and $v_{\rm 2}$, respectively. 
\item Two of the edges are denoted by $e_{\rm l}$ and $e_{\rm r}$, respectively. $e_{\rm l}$ connects $v_{\rm l}$ and $v_1$. $e_{\rm r}$ connects $v_{\rm r}$ and $v_2$.
\item Each of the remaining $b$ edges is denoted by $e_{j}$ where $1 \leq j \leq b$ is an integer. These
edges connect $v_1$ and $v_2$.
\end{itemize}
	
Let $F_{\rm l}$ and $F_{\rm r}$ be the unit spheres $S^{m-1}$. Let $\{F_j\}_{j=1}^{b}$ be a family of smooth manifolds each of which is the unit sphere or diffeomorphic to a manifold represented as a connected sum of finitely many manifolds diffeomorphic to
 the products $S^j \times S^{m-j-1}$ with integers $1 \leq j \leq m-2${\rm :} in the case $m=2$, $F_j$ must be the circle $S^1$. The connected sum is taken in the smooth category.

Then there exist an $m$-dimensional non-singular real algebraic closed and connected manifold $M$, a smooth real algebraic function $f:M \rightarrow {\mathbb{R}}$ and an isomorphism $\phi:G \rightarrow W_f$ of the graphs enjoying the following properties.
\begin{enumerate}
\item For the edge $e_{\rm l}$ of the graph $G$ and each point $p_{e_{\rm l}}$ in the interior of the edge $\phi(e_{e_{\rm l}}) \subset W_f$, the preimage ${q_f}^{-1}(p_{e_{\rm l}})$ is diffeomorphic to $F_{\rm l}$.
\item For the edge $e_{\rm r}$ of the graph $G$ and each point $p_{e_{\rm r}}$ in the interior of the edge $\phi(e_{e_{\rm r}}) \subset W_f$, the preimage ${q_f}^{-1}(p_{e_{\rm r}})$ is diffeomorphic to $F_{\rm r}$.
\item For the edge $e_{j}$ of the graph $G$ and each point $p_{e_{j}}$ in the interior of the edge $\phi(e_{j}) \subset W_f$, the preimage ${q_f}^{-1}(p_{e_j})$ is diffeomorphic to $F_j$ for each integer $1 \leq j \leq b$.
\end{enumerate}

\end{MainThm}

\begin{proof}
It is almost same as the proof of Main Theorem \ref{mthm:1}.

The only one different part is that we argue inductively and similarly starting from a connected open subset $D_0 \subset {\mathbb{R}}^2$ whose closure is a disk centered at the origin $0$. We choose the disk whose radius is $R>0$. We can choose $R$ as a sufficiently large number.

As before, we choose two distinct real numbers $t_2<t_3$ in $(-R,R)$. Put $t_1:=-R$ and $t_4:=R$.
We define the number $t_{0,j}:=\frac{t_j+t_{j+1}}{2}$ similarly only for $j=2$. 

In the first step, we consider an additional argument. 
We remove the interiors of exactly $b-1>0$ disks centered at points of the form $(t_{0,2},t_{0,2,j^{\prime}}) \in D_0$. We can choose the disks in the open set $D_0 \subset {\mathbb{R}}^2$ disjointly and as ones whose radii are all $\frac{t_3-t_2}{2}$. Note that $R$ can be chosen as a sufficiently large number. We choose a sufficiently large $R$ beforehand. 
 
Except this, we prove our statement in a similar inductive way.

Respecting this exposition completes the proof.

\end{proof}

By the construction and the structures of the maps and the manifolds, we also have the following. 
\begin{MainThm}
\label{mthm:3}
In Main Theorems \ref{mthm:1} and \ref{mthm:2}, we can construct the non-singular real algebraic closed and connected manifolds $M$ and the functions $f:M \rightarrow \mathbb{R}$ enjoying the following properties.

\begin{enumerate}
\item The manifolds $M$ are diffeomorphic to the manifolds in Examples \ref{ex:1}, \ref{ex:2} and \ref{ex:3}.
\item The functions $f$ are
represented as the compositions
of smooth real algebraic maps
which are also product-organized special generic maps into ${\mathbb{R}}^m$ with the canonical projection ${\pi}_{m,1}$.
\item The manifolds $M$ are also the zero sets of some real polynomials $F_M(x_1, \cdots x_{m+1})$.
\end{enumerate}
\end{MainThm}

\begin{figure}
	
	\includegraphics[height=75mm, width=100mm]{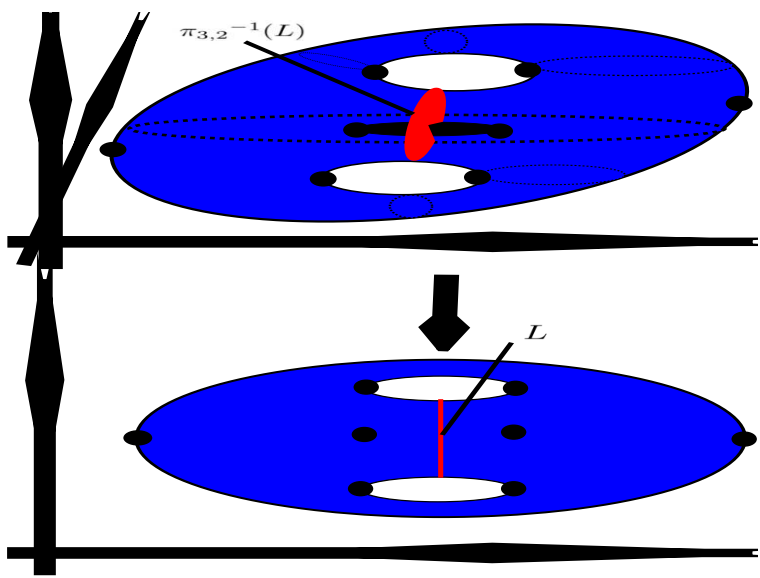}

	\caption{An example for Main Theorem \ref{mthm:2} with $m=3$ and $b=3$. These pictures show $D \subset \overline{D} \subset {\mathbb{R}}^n={\mathbb{R}}^2$ and $D^{\prime} \subset \overline{D^{\prime}} \subset {\mathbb{R}}^{n^{\prime}}={\mathbb{R}}^3$ where Theorem \ref{thm:1} is applied to show Main Theorem \ref{thm:2}. Dots are for singular points and singular values of our desired Morse function $f:M \rightarrow \mathbb{R}$. For example, the preimage ${{\pi}_{3,2}}^{-1}(L)$ of the red segment $L$ in the image of the map into ${\mathbb{R}}^2$ is diffeomorphic to $S^1 \times D^1$.}
	\label{fig:2}
\end{figure}

\begin{figure}
	
	\includegraphics[height=75mm, width=100mm]{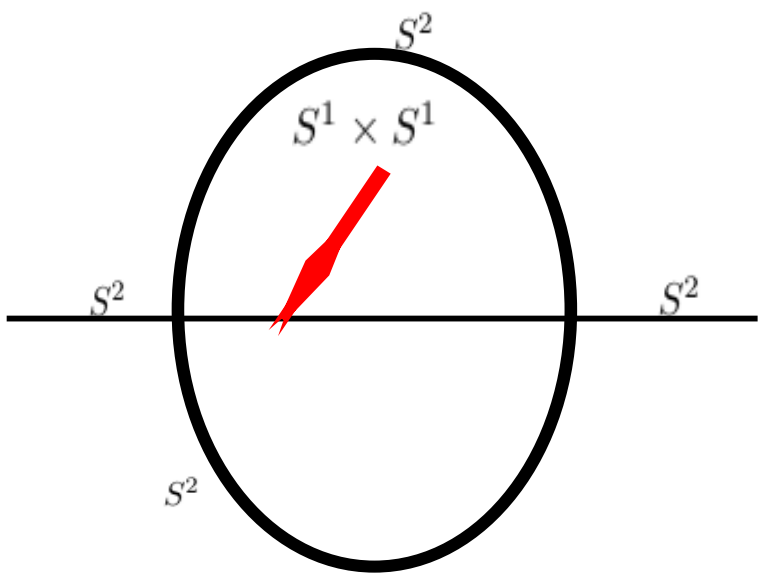}

\caption{The Reeb graph $W_f$ of the Morse function $f$ in FIGURE \ref{fig:2} and preimages. The preimage ${q_f}^{-1}(p_e)$ is diffeomorphic to $S^2$ or $S^1 \times S^1$ for a point $p_e$ in the interior of each edge $e$ of $W_f$.}
	\label{fig:3}
\end{figure}

We omit examples for Main Theorem \ref{mthm:1}. 
We present a related example for Main Theorem \ref{mthm:2} with $m=3$ and $l=3$ by FIGURES \ref{fig:2} and \ref{fig:3}. 
Here we construct a product-organized special generic map on a manifold $M$ diffeomorphic to one represented as a connected sum of exactly three copies of $S^1 \times S^2$ into ${\mathbb{R}}^3$ according to our proof. 
Here $D^{\prime}$ is diffeomorphic to a manifold represented as a boundary connected sum of exactly two copies of $S^1 \times D^2$ and one copy of $D^2 \times S^1$.
The resulting function is represented as the composition of the special generic map with the canonical projection ${\pi}_{3,1}$. 

\subsection{Previously obtained answers to Problem \ref{prob:2} and Problem \ref{prob:3}.}
We present some previously obtained answers to Problem \ref{prob:2} of the author. The following is one of results closely related to our new results.

\begin{Thm}[\cite{kitazawa1}]
\label{thm:3}
Let a graph $G$ whose vertex set and edge set are finite be given. Let $l$ be a map from the edge set into the set of all non-negative integers.
Then there exist a $3$-dimensional closed and orientable manifold $M$ and a smooth function $f:M \rightarrow \mathbb{R}$ enjoying the following properties.

	\begin{enumerate}
		\item There exists an isomorphism $\phi:G \rightarrow W_f$ from $G$ to the Reeb graph $W_f$. 
		\item For each edge $e$ of $G$, consider the element $\phi(e)$ of the Reeb graph of $W_f$ and choose an arbitrary point $p_e$ in the interior of $\phi(e)$. The preimage ${q_f}^{-1}(p_e)$ is a closed, connected and orientable surface of genus $l(e)$.
		\item If $f$ does not have a local maximum or local minimum at a singular point of $f$, then it is represented as a Morse function there for some local coordinates. 
		\item If $f$ has a local maximum or local minimum at a singular point of $f$, then there it is locally represented as either of the following four forms for suitable local coordinates.
		\begin{enumerate}
			\item A Morse function. This occurs if and only if for the singular point $p_s$, $q_f(p_s)$ is a vertex contained in the exactly one edge $e_s$ of the graph with the condition $l(e_s)=0$. More precisely, the function is regarded as a height function on a unit disk.
			\item A Morse-Bott function which is not Morse.
                        \item A function which is not Morse-Bott and which is represented as the composition of a Morse function onto the $1$-dimensional unit disk $D^1$ with a height function.
			\item A function which is not as the presented previous three cases and which is represented as the composition of a fold map onto the interior of the $2$-dimensional unit disk $D^2$ with {\rm (}the restriction of{\rm )} a height function. This occurs if and only if
 for the singular point $p_s$ of the function $f$, $q_f(p_s)$ is a vertex contained in the exactly one edge $e_s$ of the graph with the condition $l(e_s) \neq 0,1$.
	\end{enumerate}
\end{enumerate}

\end{Thm}

We present another previous result in a weaker form. This is also closely related to our new result.  

\begin{Thm}[\cite{kitazawa6}]
\label{thm:4}
Let $m \geq 2$ be an integer.
Suppose that a graph $G$ whose vertex set and whose edge set of $G$ are finite is given. Suppose that the following two maps are given.
\begin{itemize}
	\item A continuous function $g:G \rightarrow \mathbb{R}$ enjoying the following properties.
	\begin{itemize}
		\item The restriction of $g$ to the closure of each edge of the graph is injective.
		\item The function $g$ has a local maximum or a local minimum only at a vertex contained in exactly one edge in the graph $G$.
		\end{itemize}
	\item A map $l$ on the edge set of $G$ into the set of all diffeomorphism types of smooth manifolds in the following.
	\begin{itemize}
		\item The {\rm (}$m-1${\rm )}-dimensional unit sphere.
		\item Manifolds represented as connected sums of finitely many manifolds of the forms $S^{j} \times S^{m-j-1}$ with $1 \leq j \leq m-2$. Connected sums are taken in the smooth category.
	
	\end{itemize} 
Furthermore, $l(e_v)$ is the diffeomorphism type of the unit sphere if the edge $e_v$ contains a vertex $v$ where $g$ has a local maximum or a local minimum $g(v)$.
\end{itemize}
Then there exist an $m$-dimensional closed and orientable manifold $M$ and a Morse function $f:M \rightarrow \mathbb{R}$ enjoying the following properties.
\begin{enumerate}
	\item There exists an isomorphism $\phi:G \rightarrow W_f$ from $G$ to the Reeb graph $W_f$. 
	\item For each element $e$ of the edge set of $G$ and the element $\phi(e)$ of the edge set of $W_f$ and an arbitrary point $p_e$ in the interior of the edge $\phi(e)$, the diffeomorphism type of the preimage ${q_f}^{-1}(p_e)$ is $l(e)$. 
        \item It holds that $\bar{f}(\phi(v))=g(v)$ for each vertex $v$ of $G$.
\end{enumerate}
\end{Thm}
Compare these theorems to our main result. These studies are of the smooth category. We do not discuss in the real algebraic or the real analytic category in these studies.
For example, \cite{saeki2} constructs global smooth functions by using so-called {\it bump functions} or ones closely related to them. These functions are not real analytic. 
\cite{kitazawa2} also use smooth functions which are not real analytic. 

We also present our previous result on construction of smooth real algebraic functions whose Reeb graphs are isomorphic to given finite graphs. 
This is an answer to Problem \ref{prob:3} respecting Problem \ref{prob:1}.
We use terminologies such as algebraic domains and Poincar\'e-Reeb graphs of them, which are presented roughly before.
\begin{Thm}[\cite{kitazawa3}]
	\label{thm:5}
	Let $G$ be a finite graph which is also a Poincar\'e-Reeb graph of a bounded and connected algebraic domain $D \subset {\mathbb{R}}^n$ in the real affine space ${\mathbb{R}}^n$.
	Let the boundary $\partial \overline{D}$ of the closure $\overline{D}$ be a disjoint union of finitely many non-singular real algebraic hypersurfaces in ${\mathbb{R}}^n$ which are also the zero sets of some real polynomials. Let $\{S_j\}_{j=1}^l$ denote the family of the zero sets {\rm :} suppose that the family consists of exactly $l$ hypersurfaces. Let $S_j$ be the zero set of a real polynomial $f_j(x_1,\cdots x_n)$.
	Suppose also that the relations $D={\bigcap}_{j=1}^l \{x \in {\mathbb{R}}^n \mid f_j(x)> 0\}$ and $\overline{D}={\bigcap}_{j=1}^l \{x \in {\mathbb{R}}^n \mid f_j(x) \geq 0\}$ hold. Then for any integer $m \geq n$, there exist an $m$-dimensional non-singular real algebraic closed and connected manifold $M$ and a smooth real algebraic function $f:M \rightarrow \mathbb{R}$ whose Reeb graph $W_f$ is isomorphic to $G$.
	  
	\end{Thm}
	
	We can regard that the situation of Theorem \ref{thm:5} also gives a situation of Theorem \ref{thm:1}. Here we use the same notation. Furthermore, the function $f$ in Theorem \ref{thm:5} can be constructed as the composition of a smooth real algebraic map into ${\mathbb{R}}^n$ obtained by Theorem \ref{thm:1} with the canonical projection ${\pi}_{n,1}$. Here $m=n^{\prime}-1$ where $n^{\prime}$ is from Theorem \ref{thm:1}.
\begin{Ex}
	\begin{enumerate}
		\item The interior ${\rm Int}\ D^n \subset D^n \subset {\mathbb{R}}^n$ of the unit disk $D^n$ is regarded as an algebraic domain $D$ in ${\mathbb{R}}^n$ and a graph with exactly one edge and two vertices is regarded as a Poincar\'e-Reeb graph of it. The graph can be embedded naturally in ${\mathbb{R}}^n$. We also have a case for Theorem \ref{thm:5}. The function $f$ can be obtained as the canonical projection of the unit sphere $S^{m-1}$.
		\item In our proof of Main Theorem \ref{mthm:2}, a connected open set $D \subset {\mathbb{R}}^2$ obtained by removing the disjoint $b-1$ disks from the connected open set $D_0 \subset {\mathbb{R}}^2$ gives an algebraic domain of ${\mathbb{R}}^2$. The given graph $G$ in Main Theorem \ref{mthm:2} is regarded as a Poincar\'e-Reeb graph of it and we can embed this naturally in ${\mathbb{R}}^2$. We also have a case for Theorem \ref{thm:5}. The function $f$ can be also obtained as a smooth function for the case $F_j=S^{m-1}$ ($1 \leq j \leq b$) in Main Theorem \ref{mthm:2}.  

	\end{enumerate}
\end{Ex}

Last, we present a short comment. This is mainly on Main Theorem \ref{mthm:3}. It is also one on Theorem \ref{thm:4} and remaining Main Theorems. We present some in \cite{kitazawa8}.
\begin{Rem}
	\label{rem:2}
	For example, in Theorem \ref{thm:4}, consider a case such that the degree of each vertex is $1$ or $3$ and that the diffeomorphism types are those of unit spheres. In this case, we can have a Morse function $f$ such that at distinct singular points of $f$ the (singular) values are always distinct. We choose a smooth embedding into a sufficiently high dimensional Euclidean space $e:M \rightarrow {\mathbb{R}}^{n_0}$ and the map $(e,f):M \rightarrow {\mathbb{R}}^{n_0} \times \mathbb{R}={\mathbb{R}}^{n_0+1}$, which is regarded as a smooth embedding. By the celebrated theory of Nash and Tognoli, presented as "Theorem 3 of \cite{kollar}" for example, together with some theory from singularity theory of smooth maps for example, we can smoothly isotope it to a smooth real algebraic embedding by some small perturbation. By composing the canonical projection to the second $1$-dimensional real affine space $\mathbb{R}$, we also have a Morse function regarded as one obtained by considering some smooth isotopy to the original Morse function. This is regarded as a smooth real algebraic function.
For more general approximation theory in real algebraic geometry, see also \cite{bochnakcosteroy, kucharz}
	for example.
	
	We cannot apply such arguments in general. For example, in Theorem \ref{thm:4}, we consider
	Morse functions such that at some
	distinct singular points the values are same.
	We cannot apply the previous theory.
	
	Instead, we have a new method and have a smooth real algebraic function with a nice representation as in Main Theorem \ref{mthm:3}. We avoid existence
	theory and approximation theory. We do similarly in \cite{kitazawa3}.
	
	For construction of Morse functions here in the smooth category, see also \cite{michalak}. Especially, this has motivated us to present \cite{kitazawa1}.
\end{Rem}

\section{Acknowledgment.}
The author would like to thank Osamu Saeki and Shuntaro Sakurai for exciting discussions on Sakurai's master thesis \cite{sakurai}, \cite{saekisakurai}, and a talk related to them in the conference \\  https://sites.google.com/view/suzukimasahiko70/home. \\
This conference celebrates Masahiko Suzuki's 70th birthday and we would like to celebrate him again. The talk is given by Osamu Saeki. 
Their study gives an explicit smooth map on the $3$-dimensional real projective space into ${\mathbb{R}}^2$ as the restriction of a suitably chosen complex linear function on the $3$-dimensional complex space to the intersection of the $5$-dimensional unit sphere and the zero set of the complex polynomial $f(z_1,z_2,z_3)={z_1}^2+{z_2}^2+{z_3}^2$. This comes from Milnor's celebrated, interesting and well-known theory \cite{milnor}. 
The more we discuss the problems, the more we know. This leads us to \cite{kitazawa3}, followed by the present paper. The author would like to thank Osamu Saeki again for informal discussions on \cite{saeki2} with \cite{kitazawa1}. This has also motivated the author to continue related studies.

 The author would like to thank anonymous referees a lot for important comments. They have improved the paper and us. 


\begin{thebibliography}{25}
	\bibitem{bochnakcosteroy} J. Bochnak, M. Coste and M.-F. Roy, \textsl{Real algebraic geometry}, Ergebnisse der Mathematik und ihrer Grenzgebiete (3) [Results in Mathematics and Related Areas (3)], vol. 36, Springer-Verlag, Berlin, 1998. Translated from the 1987 French original; Revised by the authors.
\bibitem{bodinpopescupampusorea} A. Bodin, P. Popescu-Pampu and M. S. Sorea, \textsl{Poincar\'e-Reeb graphs of real algebraic domains}, arXiv:2207.06871.
\bibitem{bott} R. Bott, \textsl{Nondegenerate critical manifolds}, Ann. of Math. 60 (1954), 248--261.

	
\bibitem{golubitskyguillemin} M. Golubitsky and V. Guillemin, \textsl{Stable Mappings and Their Singularities}, Graduate Texts in Mathematics (14), Springer-Verlag(1974).
\bibitem{kitazawa1} N. Kitazawa, \textsl{On Reeb graphs induced from smooth functions on $3$-dimensional closed orientable manifolds with finitely many singular values}, Topol. Methods in Nonlinear Anal. Vol. 59 No. 2B, 897--912, arXiv:1902.08841.
\bibitem{kitazawa2} N. Kitazawa, \textsl{On Reeb graphs induced from smooth functions on closed or open surfaces}, Methods of Functional Analysis and Topology Vol. 28 No. 2 (2022), 127--143, arXiv:1908.04340.
\bibitem{kitazawa3} N. Kitazawa, \textsl{Real algebraic functions on closed manifolds whose Reeb graphs are given graphs}, Methods of Functional Analysis and Topology Vol. 28 No. 4 (2022), 302--308, arXiv:2302.02339, 2023. 

\bibitem{kitazawa4} N. Kitazawa, \textsl{Explicit construction of explicit real algebraic functions and real algebraic manifolds via Reeb graphs}, this is the abstract of our talk in an international conference "Algebraic and geometric methods of analysis 2023" (https://www.imath.kiev.ua/$\sim$topology/conf/agma2023/) and after a review process accepted for publication in the book of abstracts, https://imath.kiev.ua/$\sim$topology/conf/agma2023/contents/abstracts/texts/kitazawa/kitazawa.pdf.

\bibitem{kitazawa5} N. Kitazawa, \textsl{Notes on explicit special generic maps into Euclidean spaces whose dimensions are greater than $4$}, a revised version is submitted based on positive comments (major revision) by referees and editors after the first submission to a refereed journal, arxiv:2010.10078.
\bibitem{kitazawa6} N. Kitazawa, \textsl{On Reeb graphs induced from smooth functions on $3$-dimensional closed manifolds which may not be orientable}, a revised version is submitted to a refereed journal after based on positive comments by editors and referees after the second submission to a refreed journal, arXiv:2108.01300.
\bibitem{kitazawa7} N. Kitazawa, \textsl{Realization problems of graphs as Reeb graphs of Morse functions with prescribed preimages}, a revised version will be submitted to a refereed journal based on positive comments on referees and editors, arXiv:2108.06913.
\bibitem{kitazawa8} N. Kitazawa, \textsl{Some remarks on real algebraic maps which are topologically special generic maps}, submitted to a refereed journal, arXiv:2312.10646. 
\bibitem{kitazawa9} N. Kitazawa, \textsl{A note on cohomological structures of special generic maps}, a revised version is submitted based on positive comments by referees and editors after the second submission to a refereed journal.
\bibitem{kollar} J. Koll\'ar, \textsl{Nash's work in algebraic geometry}, Bulletin (New Series) of the American Mathematical Society (2) 54, 2017, 307--324.
\bibitem{kucharz} W. Kucharz, \textsl{Some open questions in real algebraic geometry}, Proyecciones Journal of Mathematics, Vol. 41 No. 2 (2022), Universidad Cat\'olica del Norte Antofagasta, Chile, 437--448.
\bibitem{martinezalfaromezasarmientooliveira} J. Martinez-Alfaro, I. S. Meza-Sarmiento and R. Oliveira, \textsl{Topological  classification of simple Morse Bott functions on surfaces}, Contemp. Math. 675 (2016), 165--179.%
\bibitem{masumotosaeki} Y. Masumoto and O. Saeki, \textsl{A smooth function on a manifold with given Reeb graph}, Kyushu J. Math. 65 (2011), 75--84.
\bibitem{maciasvirgospereirasaez} E. Mac\'ias-Virg\'os and M. J. Pereira-S\'aez, Height functions on compact symmetric spaces, Monatshefte f\"ur Mathematik 177 (2015), 119--140. 
\bibitem{michalak} L. P. Michalak, \textsl{Realization of a graph as the Reeb graph of a Morse function on a manifold}. Topol. Methods in Nonlinear Anal. 52 (2) (2018), 749--762, arXiv:1805.06727.
\bibitem{milnor} J. Milnor, \textsl{Singular points of complex hypersurfacs}, Annals of Mathematics Studies, No. 61, Princeton University Press, Princeton, N. J.; University of Tokyo Press, Tokyo, 1968.
\bibitem{moise} E. E. Moise, \textsl{Affine Structures in $3$-Manifold{\rm :} V. The Triangulation Theorem and Hauptvermutung}, Ann. of Math., Second Series, Vol. 56, No. 1 (1952), 96--114.
\bibitem{nash} J. Nash, \textsl{Real algebraic manifolds}, Ann. of Math. (2) 56 (1952), 405--421.
\bibitem{ramanujam} S. Ramanujam, \textsl{Morse theory of certain symmetric spaces}, J. Diff. Geom. 3 (1969), 213--229.
\bibitem{reeb} G. Reeb, \textsl{Sur les points singuliers d\'{}une forme de Pfaff compl\'{e}tement int\`{e}grable ou d\'{}une fonction num\'{e}rique}, Comptes Rendus
 Hebdomadaires des S\'{e}ances de I\'{}Acad\'{e}mie des Sciences 222 (1946), 847--849.
\bibitem{saeki1} O. Saeki, \textsl{Topology of special generic maps of manifolds into Euclidean spaces}, Topology Appl. 49 (1993), 265--293.
\bibitem{saeki2} O. Saeki, \textsl{Reeb spaces of smooth functions on manifolds}, International Mathematics Research Notices, maa301, Volume 2022, Issue 11, June 2022, 8740--8768, https://doi.org/10.1093/imrn/maa301, arXiv:2006.01689.
\bibitem{saekisakurai} O. Saeki and S. Sakurai, \textsl{Differentiable maps on links of complex isolated hypersurface singularities}, arXiv:2402.02365.
\bibitem{sakurai} S. Sakurai, \textsl{Singular set of the restriction of a differentiable map on Euclidean space to a submanifold}, Master's Thesis, Kyushu. Univ., 2023.
\bibitem{sharko} V. Sharko, \textsl{About Kronrod-Reeb graph of a function on a manifold}, Methods of Functional Analysis and
 Topology 12 (2006), 389--396.
\bibitem{takeuchi} M. Takeuchi, \textsl{Nice functions on symmetric spaces}, Osaka. J. Mat. (2) Vol. 6 (1969), 283--289.
\bibitem{tognoli} A. Tognoli, \textsl{Su una congettura di Nash}, Ann. Scuola Norm. Sup. Pisa (3) 27 (1973), 167--185.
\end{thebibliography}
\end{document}